\documentclass[journal]{IEEEtran}
\ifCLASSINFOpdf
\else
\fi
\usepackage{amsmath,amssymb,amsthm,bm,bbm}
\usepackage{cite}
\usepackage{amsmath}
\usepackage{xcolor}
\usepackage{graphicx}
\usepackage{flushend}

\begin{document}
%
% paper title
% Titles are generally capitalized except for words such as a, an, and, as,
% at, but, by, for, in, nor, of, on, or, the, to and up, which are usually
% not capitalized unless they are the first or last word of the title.
% Linebreaks \\ can be used within to get better formatting as desired.
% Do not put math or special symbols in the title.
\title{Netload Range Cost Curves for a Transmission-Aware Distribution System Planning under DER Growth Uncertainty}
%\title{Transmission-Aware Distribution Planning Considering DER Growth Uncertainty}
%\title{Describing Distribution-Side Planning Flexiblity for Transmission \& Distribution Coordination}
%
%
% author names and IEEE memberships
% note positions of commas and nonbreaking spaces ( ~ ) LaTeX will not break
% a structure at a ~ so this keeps an author's name from being broken across
% two lines.
% use \thanks{} to gain access to the first footnote area
% a separate \thanks must be used for each paragraph as LaTeX2e's \thanks
% was not built to handle multiple paragraphs
%

\author{Samuel~C\'ordova,~\IEEEmembership{Member,~IEEE,}
        Alexandre~Moreira,~\IEEEmembership{Member,~IEEE,}
        Miguel~Heleno,~\IEEEmembership{Senior Member,~IEEE} % <-this % stops a space
}

\maketitle

% As a general rule, do not put math, special symbols or citations
% in the abstract or keywords.
\begin{abstract}
In the face of a substantial and uncertain growth of behind-the-meter Distributed Energy Resources (DERs), utilities and regulators are currently in the search for new network planning strategies for facilitating an efficient Transmission \& Distribution (T\&D) coordination. In this context, here we propose a novel distribution system planning methodology to facilitate coordinated planning exercises with transmission system planners through the management of long-term DER growth uncertainty and its impact on the substation netload. The proposed approach is based on the design of a transmission-aware distribution planning model embedding DER growth uncertainty, which is used to determine a ``menu" of secure distribution network upgrade options with different associated costs and peak netload guarantees observed from the transmission-side, referred here as Netload Range Cost Curves (NRCCs). NRCCs can provide a practical approach for coordinating T\&D planning exercises, as these curves can be integrated into existing transmission planning workflows, and specify a direct incentive for distribution planners to evaluate peak netload reduction alternatives in their planning process. We perform computational experiments based on a realistic distribution network that demonstrate the benefits and applicability of our proposed planning approach.
\end{abstract}

% Note that keywords are not normally used for peerreview papers.

% For peer review papers, you can put extra information on the cover
% page as needed:
% \ifCLASSOPTIONpeerreview
% \begin{center} \bfseries EDICS Category: 3-BBND \end{center}
% \fi
%
% For peerreview papers, this IEEEtran command inserts a page break and
% creates the second title. It will be ignored for other modes.
\IEEEpeerreviewmaketitle

\section{Introduction} \label{sec:intro}
%%%%%%%%%%%%%%%%%%%%%%%%%%%%%%%%%%%%%%%%%%%%
%%%%%%%%%%%%%%%%%%%%%%%%%%%%%%%%%%%%%%%%%%%%
%\subsection{Motivation} \label{sec:intro_motiv}
\IEEEPARstart{D}{istribution} network upgrades constitute one of largest sources of capital expenditures in power systems, accounting for more than 30\% of the total investments costs for investor-owned utilities in the U.S. \cite{EEI2023}. In this context, distribution utilities are required to periodically justify their proposed investment plans to regulators, since these plans directly impact the resulting end-consumer electricity rates \cite{Cooke2018}. Historically, investment plans have been justified based on expected demand growth projections, for which network reinforcement needs are identified using exhaustive powerflow simulations; however, given the increasingly large integration of behind-the-meter Distributed Energy Resources (DERs), such as Rooftop Photovoltaic (RPV) and Electric Vehicles (EVs), utilities have recently started to integrate DER growth forecasts into their planning processes \cite{PGE2021, Cooke2018}. Despite the progress made so far, necessary improvements are still needed to make the current methods able to deal with DER growth uncertainty and to take into account the impact of distribution system investments on the transmission expansion planning process.

While several reports and papers have addressed the issue of incorporating DER growth forecasts into the distribution planning process, less attention has been given to DER growth uncertainty and its impact on power system planning \cite{LBNL2016, Vahid2020}. In fact, as discussed in \cite{Heymann2019, Bystrom2022, Sun2021}, the uncertainty in DER growth can play a major role when determining possible investment options for the distribution network, since different potential DER projects with diverse sizes and connecting points can result in different distribution reinforcement needs with varying locations. Furthermore, the impacts of DER growth uncertainty can also propagate to the transmission side, as each DER growth scenario can result in a different corresponding netload profile to be observed at the substation connecting Transmission and Distribution (T\&D). From the transmission perspective, this new source of netload uncertainty may potentially lead to under- or over-investments in transmission assets \cite{Sato2020}.

Ongoing discussions on T\&D coordination indicate the possible advantages of aligning distribution and transmission systems planning processes to enable a cost-effective DER integration \cite{Cooke2018, EPRI2019, IESO2021}. However, there are limitations to a practical implementation of these:

\begin{itemize}
    \item Under existing regulatory frameworks, transmission and distribution grid planning processes occur in silos, with low and non-standarized data exchange between planning entities \cite{EPRI2019}, and with desynchronized planning cycles (3-5 years for distribution \cite{PGE2021} and 10-20 years for transmission \cite{CAISO2023}).
    \item Distribution regulation does not provide sufficient incentives for reducing the distribution system's peak netload at a planning stage \cite{Kavulla2023}. Worst-case peak load conditions across distribution systems propagate upstream and tend to produce costly transmission planning requirements. This phenomenon is aggravated with the increasing levels of netload uncertainty at the planning stage, due to DER adoption.
    \item The relatively high investment cost for Non-Wire-Alternatives (NWAs), such as Battery Energy Storage Systems (BESSs), which depend on multiple value streams to be economically viable. For this purpose, an accurate valuation of the T\&D savings resulting from integrating NWAs into distribution planning (and potential impact on transmission) is fundamental \cite{IESO2021}.
    
\end{itemize}
Therefore, there is a need for practical approaches that can induce, with minimal regulatory changes, a desired level of transmission awareness when optimizing investments for distribution grids, particularly under DER growth uncertainty. 

In this context, we propose a methodology to plan distribution system investments under long-term DER growth uncertainty while providing planners at the transmission level with Netload Range Cost Curves (NRCCs). Such curves correspond to different distribution-level investments, with their implications in terms of costs and netload uncertainty at the substation, that can be presented as options to the transmission system planning.

%%%%%%%%%%%%%%%%%%%%%%%%%%%%%%%%%%%%%%%%%%%%
\subsection{Literature Review} \label{sec:intro_review}
Several research efforts have been dedicated in the last years to integrate DER uncertainty into distribution system planning. These efforts have resulted in advanced planning models embedding relevant short-term uncertainties that stem from DERs, including those associated with renewable generation \cite{Zare2018}, EV user behavior \cite{Wang2020} and demand response availability \cite{Liang2022} among others (see \cite{Ehsan2019, Vahid2020} for a thorough literature review). However, less attention has been given to the incorporation of the long-term DER growth uncertainty into distribution system planning, with only a handful of works directly addressing this issue in the context of RPV and EV growth \cite{Heymann2019, Baringo2020, Sridhar2022}. Furthermore, existing distribution system planning models that embed DER uncertainty \cite{Zare2018, Wang2020, Liang2022, Heymann2019, Baringo2020, Sridhar2022} have neglected the impact of distribution network investments on transmission-level reinforcement needs, which leaves a research gap in this regard.

A number of centralized planning models considering T\&D coordination have been proposed in the last years. These models have integrated relevant aspects such as the uncertainty related to load growth and variable renewable generation \cite{MunozDelgado2021}, the impact of distribution investment decisions on the transmission-level energy prices \cite{ElMeligy2021}, the installation of electric vehicle charging stations \cite{MoradiSepahvand2023}, and the use of open-source data and tools \cite{Muller2019}. In general, these works are based on large-scale mathematical programs that can model the operation of assets at both T\&D levels, therefore allowing an overall assessment of the impact of distribution-side investments on the transmission level. However, as a major drawback, this centralized approach directly conflicts with the current T\&D independent planning workflows and data exchange limitation imposed by the regulatory framework. While decentralized T\&D planning models have recently been proposed to address the above practical issue through the use of iterative hierarchical optimization methods \cite{Liu2018, Nikoo2020}, they do not provide a clear design of incentives that would encourage T\&D systems to collaborate, while maintaining their independence in the planning process. Furthermore, existing T\&D planning models \cite{MunozDelgado2021, ElMeligy2021, MoradiSepahvand2023, Muller2019, Liu2018, Nikoo2020} have an additional pending research gap in terms of including long-term DER growth uncertainty in their decision-making process.

In power systems operations, different T\&D coordination schemes have also been proposed in the literature \cite{Gerard2018, Givisiez2020}. In this context, the idea of using flexibility-vs-cost curves to describe the aggregate distribution-side flexibility seen from the transmission-side has gained considerable attention in the last years \cite{Silva2018, Capitanescu2023, Churkin2023}, as this approach allows T\&D operators to interact with each other in a direct yet secure fashion. However, to the best of our knowledge, this framework has only been applied to operations so far, thus leaving an untapped potential for T\&D planning applications. 

Hence, there is a need for practical network planning strategies that create a link between distribution and transmission systems planners to deal with long-term DER growth uncertainty and its impact on grid infrastructure needs at both T\&D levels. 

%%%%%%%%%%%%%%%%%%%%%%%%%%%%%%%%%%%%%%%%%%%%
\subsection{Contributions}
To bridge the gaps discussed above, we propose a methodology for distribution system planning that considers DER growth uncertainty and provides transmission system planners with different options to control substation peak netload conditions via distribution investments. These options are presented in a form of Netload Range Cost Curves (NRCCs) that can be easily integrated into transmission planning processes. %From a theoretical perspective, these curves can be interpreted as the equivalent planning version of the cost-flexibility curves proposed for TSO-DSO operation \cite{Silva2018, Capitanescu2023, Churkin2023}. 
The contributions of this paper are threefold:
\begin{enumerate}
    \item The introduction of NRCCs as a form of expressing distribution system investment options to the transmission planners, while ensuring independence in the T\&D planning process.
    \item A novel distribution system planning methodology that provides: (i) a simple and compact representation of the distribution-side's planning flexibility, which can foster a collaboration between T\&D system expansion planners; (ii) a ``menu'' of options of distribution systems plans with associated long-term netloads that can be used as an input for transmission studies.
    \item Detailed computational experiments based on a realistic distribution network. Our results demonstrate the applicability and benefits of the proposed transmission-aware distribution planning approach in the context of uncertain netload growth and adoption of RPV.
\end{enumerate}

The rest of the paper is organized as follows: Section \ref{sec:ConventionalApproaches} describes conventional approaches to address the distribution system planning problem. Section \ref{sec:proposed_methodology} presents our proposed transmission-aware distribution planning strategy embedding DER growth uncertainty. Section \ref{sec:CompExp} reports computational experiments and results, demonstrating the relevance and advantages of our proposed distribution planning approach. Finally, we draw our main conclusions in Section \ref{sec:Conclusion}.

%%%%%%%%%%%%%%%%%%%%%%%%%%%%%%%%%%%%%%%%%%%%
\section{Conventional Cost Minimization Approaches for Distribution Planning} \label{sec:ConventionalApproaches}
%%%%%%%%%%%%%%%%%%%%%%%%%%%%%%%%%%%%%%%%%%%%
The decision making process to address distribution system planning can be modeled via conventional optimization frameworks. One alternative is the use of deterministic optimization, which would render network reinforcement plans that are determined based on a cost-minimization strategy and `expected' load and DER growth projections. Another alternative is to formulate the problem in a scenario-based manner, through which a set of possible load and DER growth scenarios are incorporated so that a secure system operation can be guaranteed for all considered scenarios.

%%%%%%%%%%%%%%%%%%%%%%%%%%%%%%%%%%%%%%%%%%%%
\subsection{Deterministic Planning Model} \label{subsec:DeterministicModel}
Under an approach based on deterministic optimization, network upgrade plans are obtained considering a single user-defined load and DER growth scenario. Candidate investment options can include the reinforcement of existing network corridors as well as the installation of new voltage regulators and/or BESSs \cite{LBNL2016, Vahid2020}. The resulting deterministic planning model includes investment-related decision variables $\bm{x} = \left[ \bm{x}^B, \bm{x}^L, \bm{x}^R \right]$, operation-related decision variables $\bm{y} = \left[ \bm{\rho}^P, \bm{\rho}^Q, \bm{p}^{B,D} \bm{p}^{B,C}, \bm{q}^B, \bm{f}^P, \bm{f}^Q, \bm{e}^B, \bm{v}^{\dagger}, \bm{\tilde{v}}^R \right]$, and sets describing the hours of the day $h \in \mathcal{H}$, representative days $d \in \mathcal{D}$, distribution lines $l \in \mathcal{L}$, candidate BESSs $b \in \mathcal{B}$, and network buses $n \in \mathcal{N}$, as described next:
\begin{subequations}
\allowdisplaybreaks
\begin{align}
    & \min_{\bm{x}, \bm{y}} \; \sum_{l \in \mathcal{L}} C^{L}_l \, x_l^L + \sum_{b \in \mathcal{B}} C^B_b \, x^{B}_b + \sum_{n \in \mathcal{N}^{CR}} C^{V}_n \, x_n^{V} \label{eq:conv_obj_func} \\
    & \text{subject to:} \nonumber \\
    & \rho^P_{h,d} + \sum_{b \in \mathcal{B}(PoI)} \left( p^{B,D}_{b,h,d} - p^{B,C}_{b,h,d} \right) = \sum_{l \in \mathcal{L}^{FR}(PoI)} f_{l,h,d}^P \notag \\
    & \hspace{5.5cm} \forall h \in \mathcal{H}, d \in \mathcal{D} \label{eq:pbal_PoI} \\
    & \rho^{Q}_{h,d} + \sum_{b \in \mathcal{B}(PoI)} q^{B}_{b,h,d} = \sum_{l \in \mathcal{L}^{FR}(PoI)} f_{l,h,d}^Q \notag \\
    & \hspace{5.5cm} \forall h \in \mathcal{H}, d \in \mathcal{D} \label{eq:qbal_PoI} \\
    & - LD_{n,h,d}^{P} + \sum_{b \in \mathcal{B}(n)} \left( p^{B,D}_{b,h,d} - p^{B,C}_{b,h,d} \right) = \sum_{l \in \mathcal{L}^{FR}(n)} f_{l,h,d}^P \notag \\
    & \hspace{1.5cm} - \sum_{l \in \mathcal{L}^{TO}(n)} f_{l,h,d}^P \quad \forall h \in \mathcal{H}, d \in \mathcal{D}, n \in \mathcal{N} \label{eq:pbal_buses} \\
    & - LD_{n,h,d}^{Q} + \sum_{b \in \mathcal{B}(n)}  q^{B}_{b,h,d} = \sum_{l \in \mathcal{L}^{FR}(n)} f_{l,h,d}^Q \notag \\
    & \hspace{1.5cm} - \sum_{l \in \mathcal{L}^{TO}(n)} f_{l,h,d}^Q \quad \forall h \in \mathcal{H}, d \in \mathcal{D}, n \in \mathcal{N} \label{eq:qbal_buses} \\
    &  \Upsilon^{P}_{j} \, f_{l,h,d}^P + \Upsilon^{Q}_{j} \, f_{l,h,d}^Q \leq x^L_l \, C_l^L \notag \\
    & \hspace{3.5cm} \forall h \in \mathcal{H}, d \in \mathcal{D}, l \in \mathcal{L}, j \in \mathcal{J} \label{eq:line_cap} \\
    & e^B_{b,h,d} = e^B_{b,h-1,d} - \frac{1}{\eta^D_b} p^{B,D}_{b,h,d} + \eta^C_b \, p^{B,C}_{b,h,d} \notag \\
    & \hspace{4.5cm} \forall h \in \mathcal{H}, d \in \mathcal{D}, b \in \mathcal{B} \label{eq:bat_energy_evol} \\
    & e^B_{b,h,d} \leq DR^B_b \; x^{B}_b \hspace{2.0cm} \forall h \in \mathcal{H}, d \in \mathcal{D}, b \in \mathcal{B} \label{eq:bat_energy_limit} \\
    & \Upsilon^{P}_{j} \left( p^{B,D}_{b,h,d} - p^{B,C}_{b,h,d} \right) + \Upsilon^{Q}_{j} \, q^{B}_{b,h,d} \leq x^{B}_b \notag \\
    & \hspace{3.5cm} \forall h \in \mathcal{H}, d \in \mathcal{D}, b \in \mathcal{B}, j \in \mathcal{J} \label{eq:bat_power_limit} \\
    & - M^L_l \left( 1 - x_l^L \right) \leq \notag \\
    & v^{\dagger}_{TO(l),h,d} - v^{\dagger}_{FR(l),h,d} + 2 \left( R_l \; f^P_{l,h,d} + X_l \; f^Q_{l,h,d} \right) \notag \\
    & \hspace{0.5cm} \leq \; M^L_l \left( 1 - x_l^L \right) \quad \forall h \in \mathcal{H}, d \in \mathcal{D}, l \in \mathcal{L}^0 \subseteq \mathcal{L} \label{eq:vdrop_noreg} \\
    & - M^L_l \left( 1 - x_l^L \right) \leq \notag \\
    & \tilde{v}^{R}_{TO(l),h,d} - v^{\dagger}_{FR(l),h,d} + 2 \left( R_l \, f^P_{l,h,d} + X_l \, f^Q_{l,h,d} \right) \notag \\
    & \hspace{0.5cm} \leq \; M^L_l \left( 1 - x_l^L \right) \quad \forall h \in \mathcal{H}, d \in \mathcal{D}, l \in \mathcal{L}^R \subseteq \mathcal{L} \label{eq:vdrop_vreg} \\
    & \underline{V}_{\, n} \leq v^{\dagger}_{n,h,d} \leq \overline{V}_n \hspace{1.5cm} \forall h \in \mathcal{H}, d \in \mathcal{D}, n \in \mathcal{N} \label{eq:volt_lim} \\
    & \underline{\Phi}_{\, n}  \, \tilde{v}^{R}_{n,h,d} \leq v^{\dagger}_{n,h,d} \leq \overline{\Phi}_n \; \tilde{v}^{R}_{n,h,d} \notag \\
    & \hspace{3.2cm} \forall h \in \mathcal{H}, d \in \mathcal{D}, n \in \mathcal{N}^{R} \subseteq \mathcal{N} \label{eq:vreg_lim} \\
    &  -M^R_n \; x^R_n \leq v^{\dagger}_{n,h,d} - \tilde{v}^{R}_{n,h,d} \leq M^R_n \; x^R_n \notag \\
    & \hspace{3.2cm} \forall h \in \mathcal{H}, d \in \mathcal{D}, n \in \mathcal{N}^{CR} \subseteq \mathcal{N} \label{eq:vreg_lim_bigM} \\
    & \sum_{l \in \mathcal{L}^C(n,n')} x_l^L = 1 \quad \forall (n,n') \text{ already connected} \label{eq:one_line_cand} \\
    & \left( \, \bm{p}^{B,D}, \bm{p}^{B,C}, \bm{e}^B, \bm{v}^{\dagger}, \bm{\tilde{v}}^R \, \right) \in \mathbb{R}^{+} \label{eq:oper_vars} \\ 
    & \bm{x}^B \in \mathbb{R}^{+}, \quad \left( \bm{x}^L, \bm{x}^R \right) \in \left\{ 0,1 \right\}, \label{eq:invest_vars}
\end{align} \label{eq:conv_plan}
\end{subequations}

\noindent where $\bm{x}^L$ and $\bm{x}^V$ are binary investment decisions for new line reinforcements and voltage regulators, with associated costs $C^L$ and $C^V$; $\bm{x}^B$ is the additional BESS capacity (in MW) to be installed, with associated costs $C^B$; $\bm{\rho}^P$ and $\bm{\rho}^Q$ are the active and reactive power drawn from the transmission system at the connecting substation, respectively; $\bm{f}_l^P$ and $\bm{f}_l^Q$ are the active and reactive lineflows; $\bm{p}^{B,D}$ and $\bm{p}^{B,C}$ are the battery active power discharge and charge; $\bm{q}^B$ is the battery reactive power injection; $\bm{LD}^P$ and $\bm{LD}^Q$ are the active and reactive bus netloads; $\bm{C}^L$ is the apparent power line capacity; $\bm{e}^B$ is the energy stored (in MWh) in the BESS; $\bm{\eta}^D$ and $\bm{\eta}^C$ are the discharging and charging efficiencies; $\bm{x}^B$ is the BESS power capacity (in MW) to be installed; $\bm{DR}^B$ is the power-to-energy duration ratio (in hours); $\bm{v}^{\dagger}$ is the bus squared voltage; $\bm{R}$ and $\bm{X}$ are the line resistance and reactance; $\bm{\tilde{v}}^{R}$ is the regulator's inner squared voltage; $\underline{\bm{V}}$ and $\overline{\bm{V}}$ are the minimum and maximum bus voltage limits; $\underline{\bm{\Phi}}$ and $\overline{\bm{\Phi}}$ are the minimum and maximum voltage regulation ratio; and $\bm{M}^L$ and $\bm{M}^R$ are sufficiently large big-M parameters.

Thus, \eqref{eq:conv_plan} is a mixed-integer linear programming model, whose objective function \eqref{eq:conv_obj_func} aims to minimize the total investments costs resulting from line reinforcements and new voltage regulators and BESSs. Active and reactive nodal power balances are enforced through \eqref{eq:pbal_PoI}-\eqref{eq:qbal_PoI} for the substation that connects the distribution system to the main transmission grid, %T\&D PoI
and via \eqref{eq:pbal_buses}-\eqref{eq:qbal_buses} for the rest of the buses. Apparent power capacity limits for lines are enforced through \eqref{eq:line_cap}, in which an inner piece-wise linear representation is used based on approximating hyperplanes $j \in \mathcal{J}$ with coefficients $\Upsilon^{P}_{j}$ and $\Upsilon^{Q}_{j}$ following \cite{Masha2018}. The operation of BESSs is represented through constraints that describe the evolution of the stored energy over time \eqref{eq:bat_energy_evol}, the maximum stored energy limit \eqref{eq:bat_energy_limit} and the apparent power capacity limit \eqref{eq:bat_power_limit}. Voltage drops resulting from line powerflows are described through \eqref{eq:vdrop_noreg}-\eqref{eq:vdrop_vreg}, where $\mathcal{L}^0$ is the subset of lines whose receiving end buses do not have (neither existing nor candidate) voltage regulation, and $\mathcal{L}^R$ is the subset of lines whose receiving end buses have (either existing or candidate) voltage regulation. Bus voltage limits and voltage regulation capabilities and are enforced through \eqref{eq:volt_lim} and \eqref{eq:vreg_lim}-\eqref{eq:vreg_lim_bigM}, respectively, where $\mathcal{N}^R$ is the subset of buses with either existing or candidate voltage regulators and $\mathcal{N}^{CR}$ is a more restrict subset, which comprises the buses with non-existing but candidate voltage regulators. Constraint \eqref{eq:one_line_cand} ensures that one of the line reinforcement options (including a 'no-upgrade' option) is chosen for each network corridor, where $\mathcal{L}^C (n, n')$ indicates the set of candidate line reinforcements for each pair of connected buses $(n, n')$. Non-negative and binary decisions variables are described by \eqref{eq:oper_vars}-\eqref{eq:invest_vars}.

%%%%%%%%%%%%%%%%%%%%%%%%%%%%%%%%%%%%%%%%%%%%
\subsection{Scenario-Based Planning Model} \label{subsec:scenario-based_planning_model}

The main drawback of the {\it deterministic planning model} is the utilization of a single deterministic growth scenario, which leaves all other possible realizations of load and DER growth neglected in the decision-making process for new investments. A conventional alternative to address this issue is to formulate the investment problem within a scenario-based framework considering a set of load and DER growth scenarios as follows:
\begin{align*}
    & \min_{\bm{x}, \hat{\bm{y}}} \; \sum_{l \in \mathcal{L}} C^{L}_l \, x_l^L + \sum_{b \in \mathcal{B}} C^B_b \, x^{B}_b + \sum_{n \in \mathcal{N}^{CR}} C^{V}_n \, x_n^{V} \\
    & \text{subject to:}  \\
    & \text{Operational constraints \eqref{eq:pbal_PoI}--\eqref{eq:vreg_lim_bigM}, \eqref{eq:oper_vars}} \quad \forall k \in \mathcal{K}  \\
    & \text{Other constraints \eqref{eq:one_line_cand}, \eqref{eq:invest_vars}} 
\end{align*}

\noindent where operational decision variables and constraints are defined for each load/DER growth scenario $k \in \mathcal{K} = \left\{1, \ldots, K \right\}$, with $\hat{\bm{y}} = \left[ \bm{y}_1, \ldots, \bm{y}_k, \ldots, \bm{y}_{K} \right]$, and bus netload $\bm{LD}^P$ in \eqref{eq:pbal_buses} now depending on scenario realization $\xi_k$, i.e., $\bm{LD}^P (\xi_k)$. Thus, the above planning model determines a set of reliable distribution network upgrade options that are secure for the set of considered load/DER growth scenarios. Note, however, that this planning approach fails to capture the impact of distribution network investments on the transmission-side, which motivates the design of a transmission-aware distribution planning strategy, as described next.

%%%%%%%%%%%%%%%%%%%%%%%%%%%%%%%%%%%%%%%%%%%%
\section{Proposed Transmission-Aware Distribution Planning Strategy} \label{sec:proposed_methodology}
%%%%%%%%%%%%%%%%%%%%%%%%%%%%%%%%%%%%%%%%%%%%
As previously discussed in Sections \ref{sec:intro} and \ref{sec:ConventionalApproaches}, one of the main disadvantages of the conventional cost minimization distribution planning approaches is that they do not simultaneously consider DER growth uncertainty and its impact on the netload ``seen" from the transmission-level, which limits their practical application in coordinated T\&D planning schemes. As an alternative to provide a step towards overcoming this limitation, here, we propose a novel distribution system planning strategy integrating these two aspects so as to facilitate the dialog between distribution and transmission system planners. First, in Section \ref{subsec:planning_model}, a novel transmission-aware distribution planning model that embeds DER growth uncertainty is designed. This model can be used to determine a range of secure and efficient network upgrade options with different associated budget and peak netload intervals. Then, based on the above planning model, in Section \ref{subsec:PNvB}, we introduce the concept of NRCCs and the methodology to compute them to describe the distribution-side planning flexibility. In addition, in Section \ref{subsec:PNvB}, we include a discussion of how these curves can be implemented in current transmission planning workflows, and other practical advantages associated to the proposed planning approach.

\subsection{Transmission-Aware Distribution Planning Model} \label{subsec:planning_model}

The proposed distribution system planning model, given an available budget, aims to optimize the investment plan so as to avoid peak netloads, under the realization of different DER load growth scenarios, that largely deviate from the peak netload that would have been a result of an optimization exercise that determines investments while considering a single expected long-term DER load growth scenario. This model is written as follows:
\begin{subequations}
\begin{align}
    & \min_{\bm{x}, \hat{\bm{y}}, \bm{\lambda}} \;  W \left[ \lambda^{D} - \Lambda^{D} \right]^{+} + \left( 1 - W \right) \left[ \lambda^{R} - \Lambda^{R} \right]^{+} \label{eq:prop_obj} \\
    & \text{subject to:} \nonumber \\
    & -\lambda^{R}  \leq \rho_{h,d,k}^P \leq \lambda^{D} \quad \forall h \in \mathcal{H}, d \in \mathcal{D}, k \in \mathcal{K} \label{eq:transfer_needs}  \\
    & \left( \bm{\lambda}^D, \bm{\lambda}^R \right) \in \mathbb{R}^{+} \label{eq:pos_transfer} \\
    & \sum_{l \in \mathcal{L}} C^{L}_l \, x_l^L + \sum_{b \in \mathcal{B}} C^B_b \, x^{B}_b + \sum_{n \in \mathcal{N}^{CR}} C^{V}_n \, x_n^{V} \leq \Gamma \label{eq:max_budget} \\
    & \text{Operational constraints \eqref{eq:pbal_PoI}--\eqref{eq:vreg_lim_bigM}, \eqref{eq:oper_vars}} \quad \forall k \in \mathcal{K}  \notag \\
    & \text{Other constraints \eqref{eq:one_line_cand}, \eqref{eq:invest_vars}}, \notag
\end{align} \label{eq:prop_plan}
\end{subequations}
where $\bm{\lambda} = \left[ \lambda^D, \lambda^R \right]$, in the optimal solution, represent maximum direct (D) and reverse (R) peak netloads among all considered scenarios, which can indicate transmission-side's reinforcement needs resulting from considering multiple possible DER growth scenarios. Moreover, $\Lambda^{D}$ and $\Lambda^{R}$ are the originally `expected' direct and reverse peak netloads resulting from a deterministic DER growth projection; $\left[ \, y \, \right]^+ = \max \{ y, 0\}$ is the positive-part function; $W \in [0,1]$ is a weighting factor for prioritizing the penalization of direct or reverse peak netloads compared to the expected scenario; and, finally, $\Gamma$ is a tunable parameter that indicates the available budget to invest in the distribution grid. 

Thus, \eqref{eq:prop_plan} is a mixed-integer linear programming model. The objective function to be minimized in \eqref{eq:prop_obj} penalizes the largest deviation in terms of peak netload among all considered DER load growth scenarios compared to the expected scenario, therefore aiming at minimizing the `unexpected' transmission-side's reinforcement needs resulting from DER growth uncertainty. Constraints \eqref{eq:transfer_needs}-\eqref{eq:pos_transfer} link the active power drawn from the transmission system and the largest values of direct and reverse power flows. Constraint \eqref{eq:max_budget} limits the total investment cost based on the available budget and the last two sets of constraints ensure a secure system operation in all considered DER growth scenarios.

The distribution planning model \eqref{eq:prop_plan} allows for the simultaneous consideration of multiple possible network upgrade options based on the budget parameter $\Gamma$. In particular, by adjusting the value of $\Gamma$, different possible network upgrade plans with corresponding associated peak netload intervals and budget requirements can be obtained. This information can then be used to build the proposed NRCCs for facilitating the dialog between distribution and transmission system planners, as described next.

%%%%%%%%%%%%%%%%%%%%%%%%%%%%%%%%%%%%%%%%%%%%
\subsection{Netload Range Cost Curves} \label{subsec:PNvB}

\begin{figure}[t]
\centering
\includegraphics[width=1\linewidth]{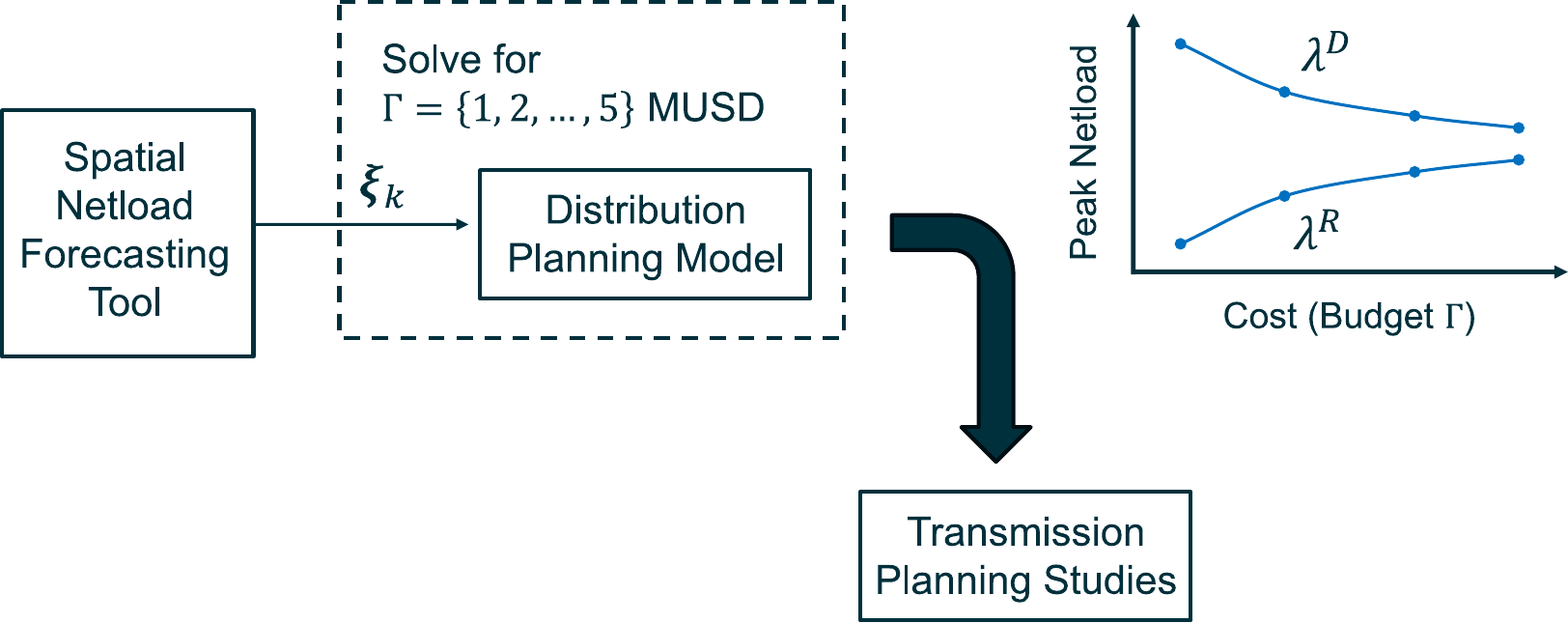}
\caption{Workflow for building NRCCs to facilitate T\&D planning coordination.}
\label{fig:PNvsB}
\end{figure}

Based on \eqref{eq:prop_plan}, NRCCs for characterizing the distribution-side planning flexibility can be built using the procedure illustrated in Fig. \ref{fig:PNvsB}. For this purpose, first, our methodology requires, as an input, a set of projected DER long-term growth scenarios $k \in \mathcal{K}$ with corresponding bus netloads $\bm{LD}^P(\xi_k)$. Although we describe a way to obtain this input in Section \ref{sec:CompExp}, proposing an ideal methodology to generate scenarios is out of the scope of this paper and our proposed approach assumes the required set of scenarios is readily available. Given the bus netload scenarios, the distribution planning model \eqref{eq:prop_plan} is solved iteratively for different budget values $\Gamma$, resulting in a range of possible direct and reverse peak netloads $\lambda^{D}$ and $\lambda^{R}$. Finally, using the above values, the NRCCs are obtained. These curves can then be shared with the transmission planning entity for its independent network studies.

\begin{figure}[t]
\centering
\includegraphics[width=1\linewidth]{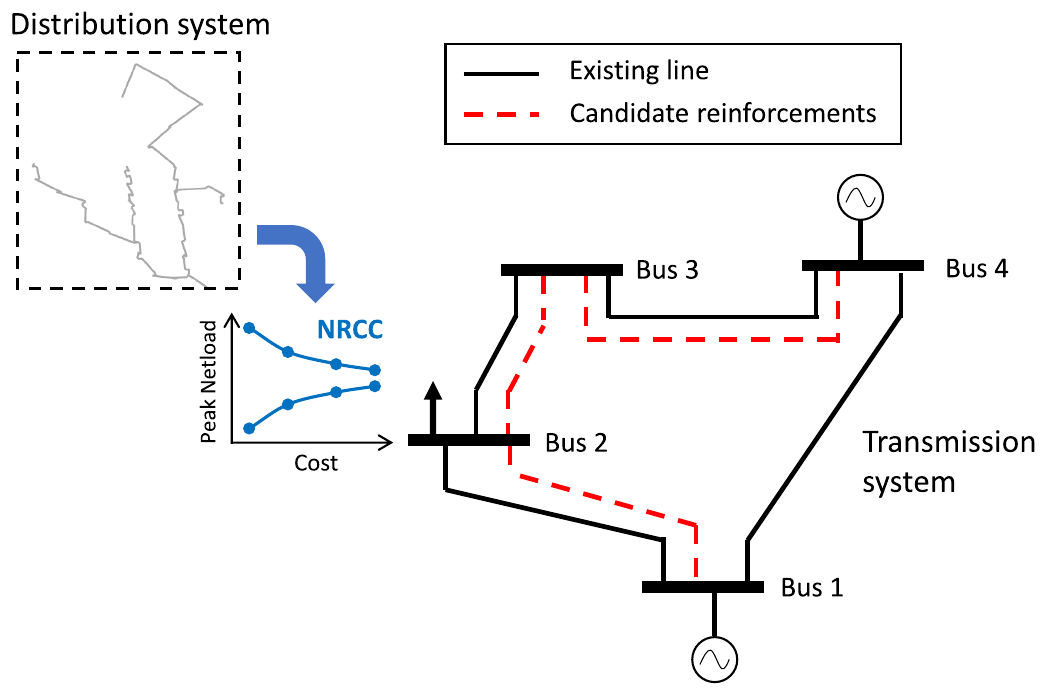}
\caption{Example of application of NRCCs in transmission planning}
\label{fig:trans_plan}
\end{figure}

One of the main advantages of the proposed NRCCs lies in their compact representation of the distribution-side's planning flexibility in a way that can be treated by the transmission counterpart. An example of the integration of NRCCs in transmission planning is illustrated in Fig. \ref{fig:trans_plan}, in which it can be observed that instead of having a single peak netload value for Bus 2, now a range of possible peak netload values for this bus is available, each of which with an associated cost. Therefore, distribution grid investments are presented as ``cost vs netload range guarantees" at the substation. 

From a transmission planner perspective, these curves describe the marginal cost of decreasing the netload range in the substation. Therefore, they can be treated equivalently to any other nodal capacity investment (such as generator or storage asset) in the expansion planning process. From a procedural perspective, the NRCCs are a non-invasive form of coordinating T\&D planning, adaptable to the two T\&D ownership models that exist in current regulatory frameworks:

\begin{enumerate}
    \item Integrated utilities: when T\&D is managed by the same entity, NRCCs can be used to communicate least-cost grid investment options from distribution to transmission planning teams.
    \item Separated utilities: when T\&D is managed by two distinct entities, NRCCs can be seen as a marginal capacity product that can be offered by the distribution company. Since marginal costs and netload range guarantees are clearly identified in NRCCs, these products are verifiable and therefore can be regulated.
\end{enumerate}

Thus, NRCCs are a straightforward and transparent mechanism for quantifying the contribution of distribution-side investments, in particular BESS and other forms of NWAs, on transmission level.

%\vspace{-0.5cm}
%%%%%%%%%%%%%%%%%%%%%%%%%%%%%%%%%%%%%%%%%%%%
\section{Computational Experiments} \label{sec:CompExp}
%%%%%%%%%%%%%%%%%%%%%%%%%%%%%%%%%%%%%%%%%%%%
In this section, we present a realistic case study to evaluate the performance of the proposed transmission-aware distribution planning approach. First, we introduce the general settings considered for the computational experiments. Then, we describe a potential method to generate realistic long-term DER growth scenarios, which results in corresponding scenarios of distribution grid nodal netloads to be considered in a planning process. It is worth emphasizing that a scenario-generation technique is not a contribution of our methodology, which assumes long-term DER growth scenarios as an available input.  Finally, we discuss the distribution planning results, starting with the network upgrade plans that result from existing distribution planning strategies, and followed by a comparison to the results obtained from implementing the proposed transmission-aware planning approach based on NRCCs.

%%%%%%%%%%%%%%%%%%%%%%%%%%%%%%%%%%%%%%%%%%%%
\subsection{General Settings} \label{sec:CompExp_settings}

\begin{figure}[t]
\centering
\includegraphics[width=0.92\linewidth]{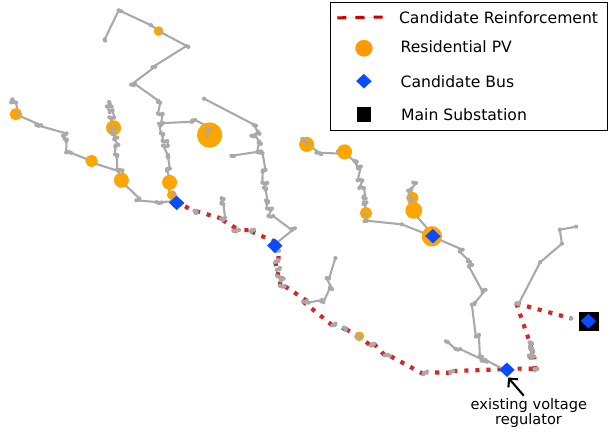}
\caption{San Francisco distribution network layout, including installed RPV capacity and candidate reinforcements. The size of the circles indicate the RPV capacity installed at each bus, ranging from 36 to 369 kW.}
\label{fig:network}
\end{figure}

Our computational experiments are based on the San Francisco distribution network provided in the SmartDS dataset \cite{SmartDS}, which includes one full year load timeseries with a 6.5 MW peak value and has a 12.47 kV base voltage, 126 existing lines, and 354 buses. This distribution system has been modified to have an initial 1.68 MW RPV capacity and an existing voltage regulator near the main substation. A 10-year planning horizon is considered in the simulations, for which three load growth scenarios with different growth rates per year are considered: Low (2\%), Medium (3\%) and High (4\%); yielding a total of 81 candidate line reinforcements, and 5 candidate buses to receive the installation of new BESS and/or voltage regulators (see Fig. \ref{fig:network}). Parameters and costs associated with candidate line reinforcements, BESSs, and voltage regulators are extracted from \cite{NREL2019, Kersting2017}.

\begin{figure}[t]
\centering
\includegraphics[width=1.01\linewidth]{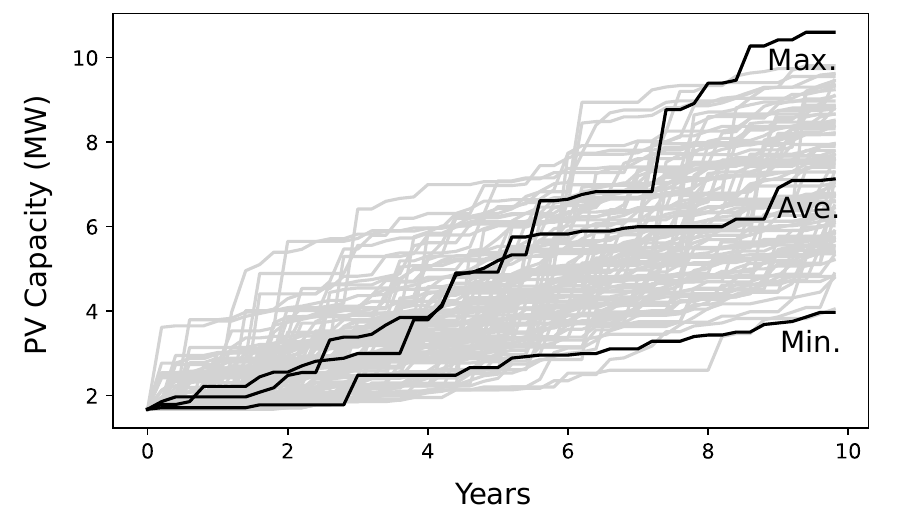}
\caption{Aggregated RPV capacity evolution resulting from 100 agent-based simulations of RPV adoption. Selected representative scenarios are highlighted in black.}
\label{fig:pv_growth}
\end{figure}

\begin{figure*}[t]
\centering
\includegraphics[width=0.95\linewidth]{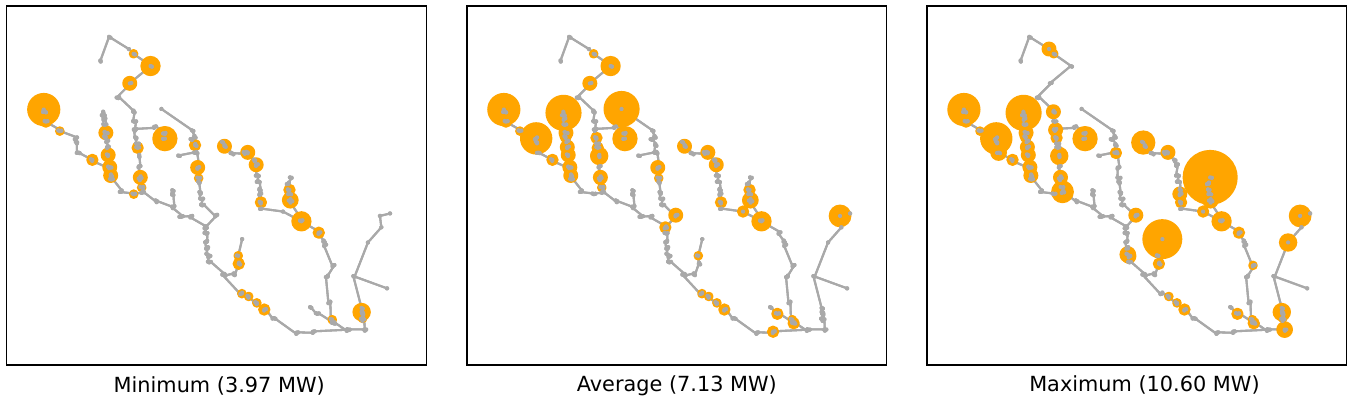}
\caption{Minimum, Average, and Maximum RPV adoption scenarios resulting from an agent-based simulation. The size of the circles indicate the installed RPV capacity at each bus, ranging from 36 to 1939 kW. The total resulting installed RPV capacity is indicated in parenthesis for each scenario.}
\label{fig:pv_scenarios}
\end{figure*}

\begin{figure}[t]
\centering
\includegraphics[width=1\linewidth]{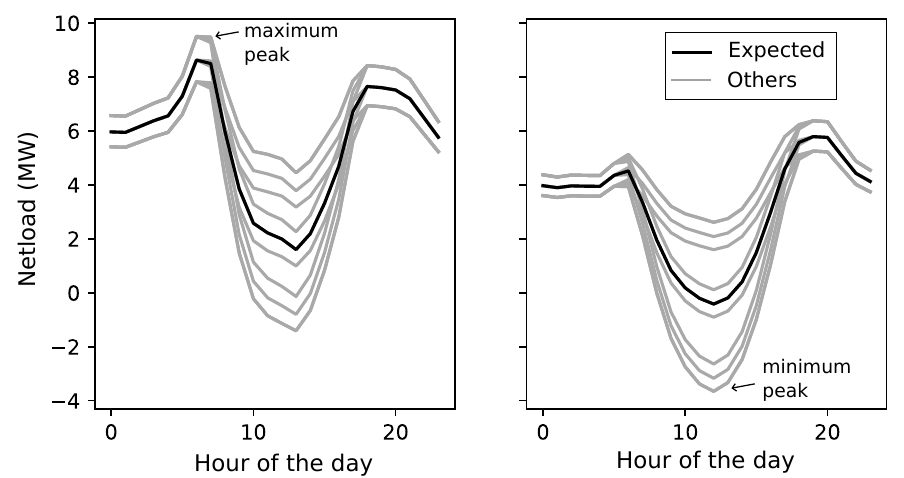}
\caption{Netload scenarios for the two chosen representative days to be used in distribution planning. Representative days were selected based on the maximum and minimum peak netloads observed across the year.}
\label{fig:repr_days}
\end{figure}

We consider a total of 100 RPV adoption scenarios to characterize realistic DER growth patterns. These scenarios are obtained from an agent-based extension of the NREL's dGen model \cite{dGen} (see Section \ref{subsec:DERscen} for more details). To reduce computational burden, 3 of the 100 DER growth scenarios are selected based on the resulting aggregate RPV capacities, yielding the Maximum, Average, and Minimum RPV adoption scenarios indicated in Figs. \ref{fig:pv_growth} and \ref{fig:pv_scenarios}. These scenarios are then combined with the 3 aforementioned load growth scenarios, therefore a total of 3x3=9 netload scenarios $\bm{\xi}_k$ are considered for this planning exercise. Such netload scenarios are further simplified by picking the two representative days that are associated with the highest and lowest peak netload values (see Fig. \ref{fig:repr_days}).

For benchmarking purposes, we compare three different distribution planning strategies in our experiments. The first one corresponds to the {\it deterministic planning model} presented in Section \ref{subsec:DeterministicModel}, which %follows current industry guidelines based on 
minimizes overall costs for distribution network upgrades while taking into account a single deterministic load/DER growth projection scenario \cite{LBNL2016, Vahid2020}. This deterministic scenario is based on the previously described Medium load growth and Average RPV adoption scenarios. The second planning strategy corresponds to the {\it scenario-based planning model} described in Section \ref{subsec:scenario-based_planning_model}, in which all 9 load/DER growth scenarios are considered in the planning process while maintaining a cost-minimization approach focused on the distribution-side. The third one is the {\it proposed transmission-aware distribution planning strategy}, described in Section \ref{sec:proposed_methodology}, in which a transmission-aware planning model embeds DER growth uncertainty and is used to build NRCCs that describe possible distribution network upgrade options with different associated budgets and peak netload intervals. Under this strategy, a symmetric weighting factor of $W=0.5$ is considered for both direct and reverse peak netloads.%, and the originally `expected' transmission-needs $\Lambda^D$ and $\Lambda^R$ are obtained using the maximum and minimum peak netloads resulting from the `Conventional' plannning strategy.

The conventional and proposed approaches are implemented in Python, with Pyomo and Gurobi building and solving the optimization problems, and our computational experiments have been performed on a PC with an Intel Core i7 2.20-GHz processor and 32 GB of RAM under a 64-bit Windows 11 operating system. RPV projects are simulated using the PySAM package \cite{PySAM}. %Optimization problems are solved using the Pyomo package \cite{Pyomo} and Gurobi solver \cite{gurobi}. 
In addition, we run AC powerflows for result validation using the OpenDSSDirect package \cite{OpenDSSDirect}.

%%%%%%%%%%%%%%%%%%%%%%%%%%%%%%%%%%%%%%%%%%%%
\subsection{DER Growth Scenario Generation} \label{subsec:DERscen}

Here we describe the methodology used to generate realistic DER growth scenarios for distribution planning. This methodology is based on simulating future RPV adoption patterns, while considering the uncertainty related to the location and size of the new RPV projects in the distribution grid.

Modern RPV adoption models are typically based on Innovation Diffusion theory, which characterizes how a new technology is adopted over time. Under this approach, RPV adoption is characterized as follows \cite{dGen, Kiesling2012}:
\begin{align}
\frac{dF(t)}{dt} \frac{1}{1-F(t)} = p + q \, F(t),
\label{eq:agg_diff_cont}
\end{align}

\noindent where $F(t)$ is the fraction of RPV adopters at time $t$, and $p$ and $q$ are the coefficient of innovation and imitation, respectively. Model \eqref{eq:agg_diff_cont} has successfully been implemented to describe aggregate RPV adoption trends in the U.S. \cite{dGen}; however, in its current version, is incapable of modeling the geographical distribution of the new RPV projects in a distribution network, as it only models aggregate dynamics. To address this issue, we modified \eqref{eq:agg_diff_cont} for its application in an agent-based simulation framework, as described next.

\begin{figure}[t]
\centering
\includegraphics[width=1.0\linewidth]{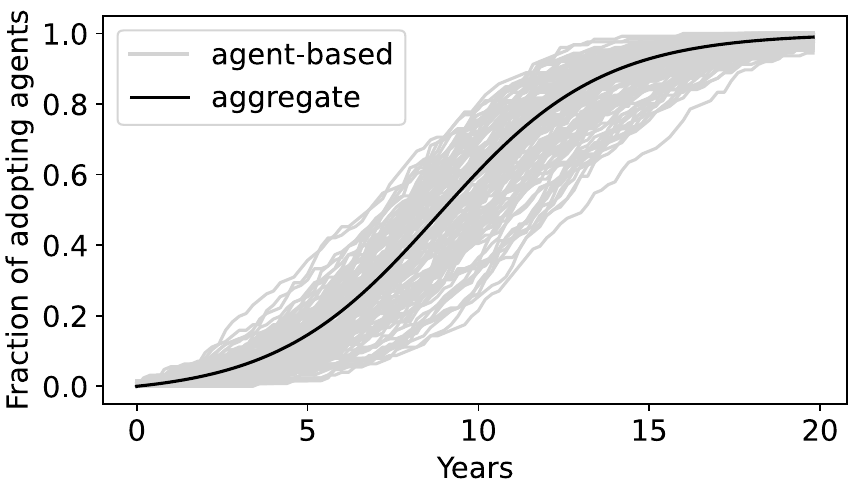}
\caption{Evolution of fraction of adopting agents over time for aggregate continous model (blue) and multiple agent-based simulations (grey).}
\label{fig:BassGrowth}
\end{figure}

As a first step, each bus $n \in \mathcal{N} = \{1, \ldots, N\}$ in the distribution network is treated as an independent agent, for which an adoption status $x_{n,t}$ at time $t$ is defined ($1$ if adopted, $0$ if not). Thus, based on \eqref{eq:agg_diff_cont} and \cite{Kiesling2012}, the probability of adoption at time $t +\Delta t$ can be computed as follows:

\begin{align}
Pr \left( x_{n,t +\Delta t} = 1 \, \lvert \, x_{n,t}=0 \right) \, = \left( p + \frac{q}{N} \sum_{n \in \mathcal{N}} x_{n,t} \right) \Delta t \label{eq:prob_adopt}
\end{align} 
which can be used to run Monte-Carlo simulations describing potential adoption scenarios for each agent $n \in \mathcal{N}$. As an example, the results of running 100 independent agent-based simulations is shown in Fig. \ref{fig:BassGrowth}, in which the fraction of adopting agents over time for each simulated scenario is presented, and compared to the one estimated by the aggregate continuous model \eqref{eq:agg_diff_cont}.

Another relevant output of the modified agent-based simulation framework is the size (in kW) of the RPV project assigned to each agent. To capture the economic rationality of the agents, this aspect is modeled by solving an optimization problem in which the RPV capacity that maximizes the Net Present Value (NPV) of each potential project is calculated independently. Thus, when a new RPV project in installed into the distribution grid (i.e., $x_{n,t} = 1$), this pre-determined RPV capacity value is used. Furthermore, similar to \cite{dGen}, a maximum market share (MMS) dependent on the project's NPV and payback period is also considered here. Within the modified agent-based framework, this means first running a Bernoulli trial (with a success probability equal to the MMS) to determine the subset of agents $n \in \mathcal{N}' \subset \mathcal{N}$ that will be considered for an eventual adoption in the future resulting from \eqref{eq:prob_adopt}. 

The  Minimum, Average, and Maximum RPV scenarios resulting from running 100 agent-based simulations and picking the 3 most representative cases (in terms of aggregate RPV capacity) are illustrated in Figs. \ref{fig:pv_growth} and \ref{fig:pv_scenarios}. Observe that not only the total installed RPV capacity varies significantly across the 3 scenarios (from 3.97 to 10.60 MW), but also the spatial distribution of the new RPV projects. Therefore there is a two-dimensional RPV growth uncertainty related to the size and location of new projects. As described in Section \ref{sec:CompExp_settings}, the above 3 RPV adoption scenarios are the ones considered for distribution planning.

%%%%%%%%%%%%%%%%%%%%%%%%%%%%%%%%%%%%%%%%%%%%
\subsection{Planning Results:Deterministic Planning Model} \label{subsec:conv_plan}

\begin{figure}[t]
\centering
\includegraphics[width=1.0\linewidth]{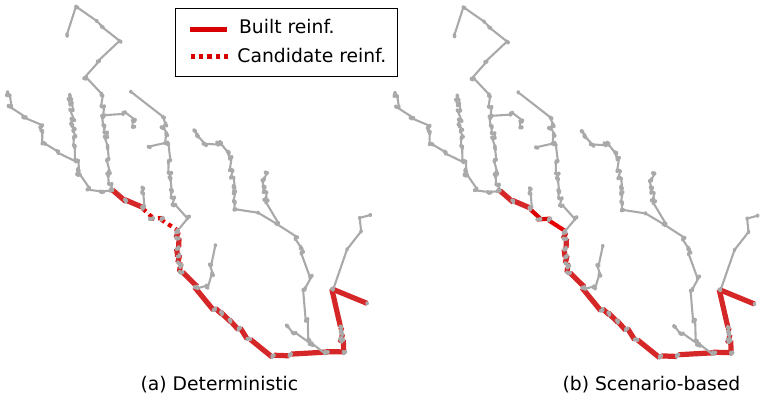}
\caption{Network upgrade plan resulting from (a) the {\it deterministic planning model} and (b) the {\it scenario-based planning model}.}
\label{fig:conv_rob_plan}
\end{figure}

The results obtained from implementing the {\it deterministic planning model} are presented in Fig. \ref{fig:conv_rob_plan}a. Under this approach, a total investment cost of 0.98 MUSD is obtained as a result of line reinforcements. Moreover, minimum and maximum peak netload values for the single considered long-term scenario are equal to -0.42 MW and +8.62 MW, respectively. These observed results stem from the {\it deterministic planning model}'s objective to minimize investment costs only. Hence, the cheapest network upgrade option is preferred (line reinforcements in this case) without considering the potential impact of the distribution system's peak netload on the transmission-side.

\begin{figure}[t]
\centering
\includegraphics[width=0.9\linewidth]{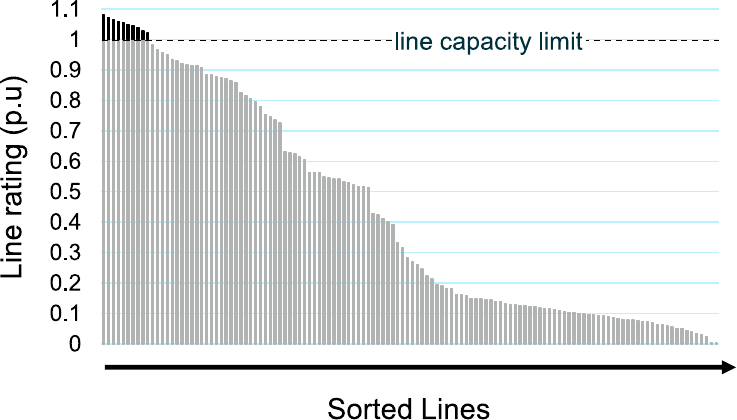}
\caption{Line rating in per unit (p.u.) observed during the maximum peak hour under a High load and DER growth scenario. Line capacity threshold is indicated to identify limit violations.}
\label{fig:line_rating}
\end{figure}

The above planning results were tested under all 9 load and DER growth scenarios by running detailed AC powerflow simulations, which revealed significant line capacity violations (see Fig. \ref{fig:line_rating}). These violations result from the fact that only a single expected scenario of load growth and RPV adoption is considered under the framework of the {\it deterministic planning model}, which underestimates the network upgrades required for most extreme scenarios with higher load and RPV growth.

%%%%%%%%%%%%%%%%%%%%%%%%%%%%%%%%%%%%%%%%%%%%
\subsection{Planning Results: Scenario-Based Planning Model} \label{subsec:stoch_plan}

The results obtained from implementing the {\it scenario-based planning model} are presented in Fig. \ref{fig:conv_rob_plan}b, for which all selected 9 load/DER growth scenarios are considered in the planning process. Observe that unlike indicated by the solution attained via the {\it deterministic planning model} (see Fig. \ref{fig:conv_rob_plan}a), additional distribution line reinforcements are required in this case to ensure a secure system operation under the realization of all selected scenarios. As a result, the total investment cost is $1.24 > 0.98$ MUSD. Moreover, note that in this case, a wider peak netload range equal to [-3.66, +9.50] MW is obtained, which is significantly larger than the one originally estimated by the {\it deterministic planning model}, i.e., [-0.42, +8.62] MW. This larger range of peak netloads is a natural consequence of considering different long-term scenarios.

%\vspace{-0.5cm}
%%%%%%%%%%%%%%%%%%%%%%%%%%%%%%%%%%%%%%%%%%%%
\subsection{Planning Results: Proposed Transmission-Aware Distribution Planning Strategy}

The performance of the {\it proposed transmission-aware distribution planning strategy} described in Section \ref{sec:proposed_methodology} is tested by solving planning model \eqref{eq:prop_plan} considering all 9 load/DER growth scenarios for different available budget values $\Gamma$. A direct and reverse `expected' transmission reinforcement need of $\Lambda^D = 8.62$ MW and $\Lambda^R = 0.42$ MW is used in \eqref{eq:prop_plan}, based on the peak netloads estimated by the {\it deterministic planning model} (see Section \ref{subsec:conv_plan}).

\begin{figure*}[t]
\centering
\includegraphics[width=1.0\linewidth]{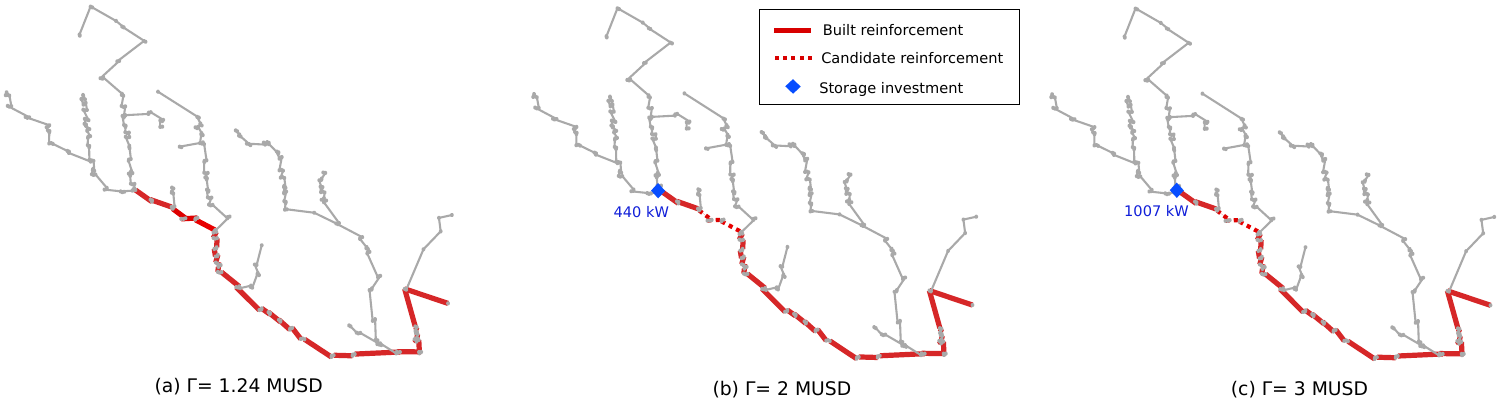}
\caption{Network upgrade plans resulting from {\it proposed transmission-aware distribution planning strategy} for different available budget $\Gamma$ values}
\label{fig:all_reinf}
\end{figure*}

As a first analysis, the planning model \eqref{eq:prop_plan} is solved considering $\Gamma =  1.24$ MUSD, which corresponds to the total investment cost obtained via the {\it scenario-based planning model}, as discussed in Section \ref{subsec:stoch_plan}. The planning results for this case are illustrated in Fig. \ref{fig:all_reinf}a, which as anticipated, are consistent with the ones obtained through the {\it scenario-based planning model} (see Fig. \ref{fig:conv_rob_plan}b), including the resulting direct and reverse peak netloads $\lambda^D = 9.50$ and $\lambda^R = 3.66$ MW. %stemming from the distribution system's peak netload values.

As a second analysis, the available budget was increased to $\Gamma = 2$ MUSD, which resulted in the network upgrade plan presented in Fig. \ref{fig:all_reinf}b. Observe how unlike the case with $\Gamma = 1.24$ MUSD (see Fig. \ref{fig:all_reinf}a), a 440kW BESS is installed in this new case, which translates into reduced direct and reverse peak netloads of $\lambda^D = 9.06 < 9.50$ MW and $\lambda^R = 3.22 < 3.66$ MW. In addition this solution indicated the deferral of some of the originally required distribution line upgrades due to the BESS smoothing dispatch. We observe similar results by further raising the available budget to $\Gamma = 3$ MUSD. As a consequence of this larger budget, the BESS capacity is increased to $1007 > 440$ kW, and the peak netloads are further reduced to $\lambda^D = 8.62 < 9.06$ MW and $\lambda^R = 2.65 < 3.22$ MW (see Fig. \ref{fig:all_reinf}c), therefore very likely decreasing the levels of needed investments in the transmission system.

\begin{figure}[t]
\centering
\includegraphics[width=1.0\linewidth]{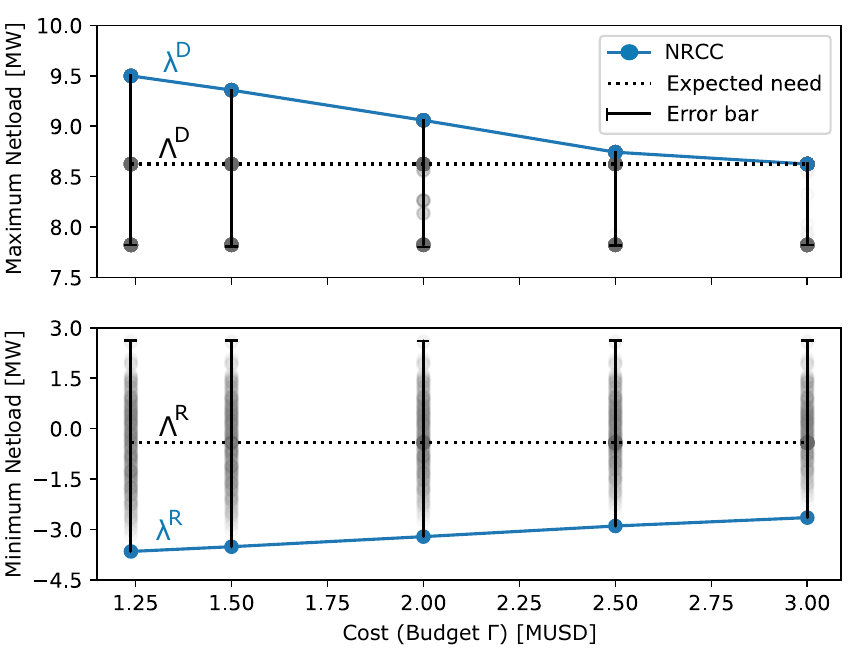}
\caption{Netload Range Cost Curves (NRCC) resulting from distribution planning for different costs (budget $\Gamma$). Dispersion bars resulting from implementing the planning solutions across all 300 DER penetration and load growth scenarios are included.}
\label{fig:PNvB_res}
\end{figure}

A summary of the resulting direct and reverse peak netloads $\lambda^D$ and $\lambda^R$ for different available budgets $\Gamma$ is presented in Fig. \ref{fig:PNvB_res}, which corresponds to the NRCCs originally described in Section \ref{subsec:PNvB}. It is worth noting how peak netloads $\lambda^D$ and $\lambda^R$ decrease as a larger budget $\Gamma$ is made available for the distribution network, ultimately converging to the `expected' requirements $\Lambda^D$ and $\Lambda^R$ originally estimated by the {\it deterministic planning model}. Such reduction results from the installation of new BESSs in the distribution network, which reduces the system's minimum and maximum peak netloads across all the selected 9 load/DER growth scenarios. 
Note that Fig. \ref{fig:PNvB_res} also includes dispersion bars to depict the peak netloads obtained by implementing and assessing the performance of the attained planning decisions (which were determined while considering 9 scenarios as previously mentioned) under 
the realization of 300 scenarios of DER penetration and load growth. These 300 scenarios are a combination of 
the original 100 RPV adoption scenarios (see Fig. \ref{fig:pv_growth}) with the 3 load growth scenarios under consideration. As can be seen, the dispersion bars indicate that the peak netload for all scenarios are within the bounds of the NRCC, which were built using only the Maximum, Average, and Minimum RPV adoption scenarios.
\color{black}

From the observation of Fig~\ref{fig:PNvB_res}, it is possible to derive the three main characteristics of the proposed NRCCs: 

\begin{enumerate}
    \item The notion of distribution-side planning flexibility, i.e. the exploration of feasible investment solutions that go beyond what is strictly necessary to ensure distribution grid security to decrease the uncertainty ``seen" from the transmission side;
    \item The expression of this planning flexibility in the form of netload range guarantees at the distribution feeder head or substation transformer;
    \item The ability to link these netload range guarantees to concrete tractable distribution investments and corresponding incremental costs, which facilitates the regulatory processes around T\&D. 
\end{enumerate}

\begin{figure}[t]
\centering
\includegraphics[width=1.0\linewidth]{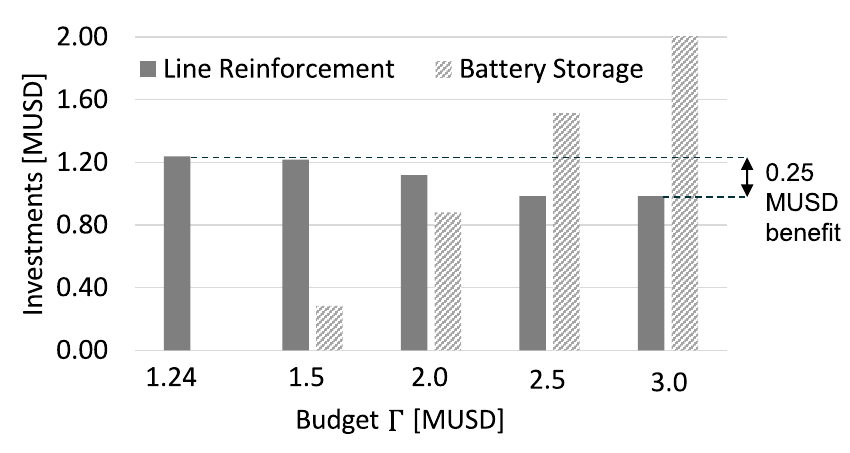}
\caption{Investment cost breakdown for different available budget $\Gamma$ values}
\label{fig:invest_bar}
\end{figure}

After disaggregating the planning solutions by investment type, it is possible to observe the benefit of strategically locating new BESS investments to avoid reconductoring some of the distribution's network corridors (see Fig. \ref{fig:all_reinf}). This benefit can directly be quantified, as illustrated in Fig. \ref{fig:invest_bar}, in which it can be observed that efficiently-located BESS investments result in a reduction of line reinforcement needs at a distribution level. In particular, for $\Gamma=3$ MUSD, line savings equal to 0.25 MUSD are obtained as a result of strategically locating a 1007kW BESS in the middle of the feeder. This illustrates how, in our proposed framework, the value of strategically locating BESS (and other possible NWAs) is captured from a T\&D point of view. 

Thus, from the case study above, it results clear that NRCCs provide a series of practical benefits that facilitate their implementation in the context of a coordinated T\&D planning. First, they represent the distribution-side's planning flexibility in a simple and compact manner, which facilitates its integration within existing transmission planning workflows. Second, from a distribution grid perspective, NRCCs are able to distinguish assets strictly required to ensure system security from those that can be installed to additionally benefit the transmission side.  

%\vspace{0cm}
%%%%%%%%%%%%%%%%%%%%%%%%%%%%%%%%%%%%%%%%%%%%
\section{Conclusion and Future Work} \label{sec:Conclusion}
%%%%%%%%%%%%%%%%%%%%%%%%%%%%%%%%%%%%%%%%%%%%

Driven by a large and uncertain growth of behind-the-meter DERs, utilities and regulators are currently in the search for new practical network planning strategies capable of guaranteeing an efficient yet reliable T\&D system expansion. Motivated by this challenge, we presented in this paper a novel distribution system planning strategy for facilitating T\&D coordination through the management of DER growth uncertainty and its impact on transmission-side reinforcement needs. The proposed strategy is based on a novel transmission-aware distribution planning model, which introduces the concept of NRCCs to describe the distribution-side planning flexibility that can result from allocating additional resources to the distribution system so as to reduce peak netload uncertainty, and potentially defer transmission-side investments. Computational experiments on a realistic distribution network demonstrated the benefits and applicability of the proposed planning approach, which would require minimal changes to the current regulatory settings. These benefits include:

\begin{itemize}
    \item A simple and compact description of distribution-side planning flexibility through NRCCs, which allow a seamless integration within current transmission planning workflows for harvesting the benefits of a coordinated T\&D planning.
    \item The provision of a direct incentive for distribution planners to evaluate peak netload reduction alternatives, in which higher budgets are allocated to the distribution-side in its role of supporting transmission planning.
    \item A straightforward and transparent mechanism to capture the benefits of distribution investments (including NWAs) that can be introduced in the current regulatory distribution system planning processes, regardless of the utility ownership mechanisms.
\end{itemize}

Potential future works could focus on a transmission planning with consideration of NRCCs curves and quantifying the monetary savings resulting from this additional source of planning flexibility.

%\vspace{-0 cm}
% if have a single appendix:
%\appendix[Proof of the Zonklar Equations]
% or
%\appendix  % for no appendix heading
% do not use \section anymore after \appendix, only \section*
% is possibly needed

% use appendices with more than one appendix
% then use \section to start each appendix
% you must declare a \section before using any
% \subsection or using \label (\appendices by itself
% starts a section numbered zero.)
%

% Can use something like this to put references on a page
% by themselves when using endfloat and the captionsoff option.
\ifCLASSOPTIONcaptionsoff
  \newpage
\fi

% trigger a \newpage just before the given reference
% number - used to balance the columns on the last page
% adjust value as needed - may need to be readjusted if
% the document is modified later
%\IEEEtriggeratref{8}
% The "triggered" command can be changed if desired:
%\IEEEtriggercmd{\enlargethispage{-5in}}

% references section

% can use a bibliography generated by BibTeX as a .bbl file
% BibTeX documentation can be easily obtained at:
% http://mirror.ctan.org/biblio/bibtex/contrib/doc/
% The IEEEtran BibTeX style support page is at:
% http://www.michaelshell.org/tex/ieeetran/bibtex/
%\bibliographystyle{IEEEtran}
% argument is your BibTeX string definitions and bibliography database(s)
%\bibliography{IEEEabrv,../bib/paper}
%
% <OR> manually copy in the resultant .bbl file
% set second argument of \begin to the number of references
% (used to reserve space for the reference number labels box)
%\vspace{-0.2cm}
\bibliographystyle{IEEEtran}
\bibliography{review}

% biography section
% 
% If you have an EPS/PDF photo (graphicx package needed) extra braces are
% needed around the contents of the optional argument to biography to prevent
% the LaTeX parser from getting confused when it sees the complicated
% \includegraphics command within an optional argument. (You could create
% your own custom macro containing the \includegraphics command to make things
% simpler here.)
%\begin{IEEEbiography}[{\includegraphics[width=1in,height=1.25in,clip,keepaspectratio]{mshell}}]{Michael Shell}
% or if you just want to reserve a space for a photo:

% You can push biographies down or up by placing
% a \vfill before or after them. The appropriate
% use of \vfill depends on what kind of text is
% on the last page and whether or not the columns
% are being equalized.

%\vfill

% Can be used to pull up biographies so that the bottom of the last one
% is flush with the other column.
%\enlargethispage{-5in}

% that's all folks
\end{document}